\documentclass[english]{amsart}

\usepackage{babel}
\usepackage{amstext}
\usepackage{amsmath}
\usepackage{amsfonts}
\usepackage{latexsym}
\usepackage{ifthen}
\usepackage{xypic}
\xyoption{all}
\renewcommand{\O}{{\mathcal O}}
\newcommand{\I}{{\mathcal I}}
\newcommand{\F}{{\mathcal F}}
\newcommand{\LB}{{\mathcal L}}
\newcommand{\Z}{{\mathcal Z}}
\newcommand{\g}{{\mathfrak g}}
\newcommand{\PN}{{\mathbb P}}

\newcommand{\codim}{{\rm codim}}
\newcommand{\Pic}{{\rm Pic}}

\newcommand{\lra}{\longrightarrow}
\newcommand{\KC}{{\mathbb C}}

\newcommand{\Sl}{{\rm Sl}}
\newcommand{\im}{{\rm im}}
\newcommand{\ad}{{\rm ad}}
\newcommand{\z}{\mathfrak z}

\newtheorem{lemma1}[equation]{}

\newenvironment{lemma}{\begin{lemma1}{\bf Lemma.}}{\end{lemma1}}

\newenvironment{abs}{\begin{lemma1}\rm}{\end{lemma1}}

\newenvironment{proposition}{\begin{lemma1}{\bf Proposition.}}{\end{lemma1}}

\newenvironment{remark}{\begin{lemma1}{\bf Remark.}\rm}{\end{lemma1}}

\begin{document}

\title{Fano threefolds with sections in $\Omega_V^1(1)$}
\author{Priska Jahnke}
\author{Ivo Radloff}
\address{Mathematisches Institut \\ Universit\"at Bayreuth \\ D--95440 Bayreuth/
Germany}
\email{priska.jahnke@uni-bayreuth.de}
\email{ivo.radloff@uni-bayreuth.de}
\date{\today}
\thanks{The authors were supported by the Deutsche Forschungsgemeinschaft.}
\maketitle

\section*{Introduction}

Let $V$ be a Fano manifold of Picard number one, and let $\O_V(1)$ be an ample generator of $\Pic(V)$. Usually $H^0(V, \Omega_V^1(1)) = H^0(V, \Omega_V^1 \otimes \O_V(1)) = 0$. The existence of a form $0 \not= \theta \in H^0(V, \Omega_V^1(1))$ is therefore a special condition. Two particular cases are well known: firstly, if  $\dim V = 2r+1$ is odd and $\theta \in H^0(V, \Omega_V^1(1))$ induces a bundle sequence
 \begin{equation} \label{contactseq}
   0 \lra F \lra T_V \stackrel{\theta}{\lra} \O_V(1) \lra 0
 \end{equation}
with maximal non--integrable kernel $F$, then $V$ is a so called Fano contact manifold, and it is conjectured that $V$ is homogeneous in this case. Secondly, if $0 \not= \theta \in H^0(V, \Omega_V^1(1))$, and $d\theta \wedge \theta \in H^0(V, \bigwedge^3\Omega_V^1 \otimes \O_V(2))$ is the zero section, then the kernel of $\theta$ induces a foliation on $V$, which is again a quite special situation.

\vspace{0.2cm}

In general, a section $\theta \in H^0(V, \Omega_V^1(1))$ will neither induce a bundle sequence like (\ref{contactseq}), nor will $d\theta \wedge \theta \in H^0(V, \bigwedge^3\Omega_V^1 \otimes \O_V(2))$, the section deciding on integrability, be either free of zeroes or completely zero. In general, both $\theta$ and $d\theta \wedge \theta$ will have nontrivial vanishing loci, and the interesting question is in how far these reflect the geometry of $V$. We note that $\theta$ cannot vanish on a divisor, since $V$ has no holomorphic $1$--forms by assumption. 

Using Iskovskikh's classification of Fano threefolds, the coarse picture is as follows:

\vspace{0.2cm}

\noindent {\bf Theorem.} {\em Let $V$ be a Fano threefold of Picard number one and index $r$, and denote by $\O_V(1)$ an ample generator of $\Pic(V)$. If we have on $V$ a holomorphic section $0 \not= \theta \in H^0(V, \Omega_V^1(1))$, then $V$ is in one of the following families
\begin{enumerate}
  \item $V_{22}$. If $V$ is general in the family, then $H^0(V, \Omega_V^1(1)) = \KC^3$ and $d\theta \wedge \theta \in H^0(V, \O_V(1))$ vanishes. 
  \item $V_{18}$. Here $H^0(V, \Omega_V^1(1)) = \KC$ and $d\theta \wedge \theta \in H^0(V, \O_V(1))$ is non--vanishing for any member of the family.
\end{enumerate}}

\vspace{0.2cm}

For particular members of the family, where a more detailed description of the Fano manifold in question is available, we can say far more. In 1.), the special member $V_{22}^s$, the Mukai--Umemura threefold, is almost homogeneous. Here $H^0(V_{22}^s, \Omega_{V_{22}^s}^1(1)) = \KC^3$, and $d\theta \wedge \theta \in H^0(V_{22}^s, \O_{V_{22}^s}(1))$ either cuts out precisely the divisor of lines on $V_{22}^s$, or vanishes completely, defining an almost homogeneous foliation. In contrast to this special case, $d\theta \wedge \theta$ always vanishes on a general $V_{22}$. 

By Mukai's classification, a $V_{18}$ is a complete intersection of two hyperplanes in a $5$ dimensional homogeneous contact manifold $M$. Here the space $H^0(V, \Omega_V^1(1)) = \KC$ is simply generated by the pull back of the contact form on $M$ to $V$. On $M$, using the contact sequence, vector fields and hyperplane sections may be identified. If we think in this way of $V$ being the complete intersection of the hyperplanes corresponding to $X_1, X_2 \in H^0(M, T_M)$, then $H^0(V, \Omega_V^1(1))$ is generated by the restriction of $X_1 \wedge X_2$, and the nonvanishing section that decides on integrability corresponds to $[X_1, X_2]$.


\section{Existence of sections in $\Omega_V^1(1)$}
\setcounter{equation}{0}

We will use both Iskovskikh's and Mukai's classification to determine all Fano threefolds $V$ with Picard number one which admit a holomorphic section in $\Omega_V^1(1)$. For the convenience of the reader we have added the classification from \cite{Isk} and \cite{Mukai} in the appendix. 

Some notations: denote the {\em index} of $V$ by $r$, i.e. $-K_V = rH$, where $\O_V(1) = \O_V(H)$ is the fundamental divisor on $V$. By Kobayashi and Ochiai's criterion, $1 \le r \le 4$ and $r = 3,4$ if and only if $V \simeq Q_3,\PN^3$, respectively. It remains hence to classify the cases $r = 1$ and $r = 2$. Let $d = H^3$ be the {\em degree} of $V$. A Fano threefold of degree $d$ and index $1$ we call $V_d$, by $V_{2,d}$ we denote a Fano threefold of index $2$ and degree $d$.

\vspace{0.2cm}

Iskovskikh uses the method of double projection from a line for his classification. The existence of lines was proved by Shokurov in \cite{Shokurov}. Key of Iskovskikh's method is \cite{Isk}, Theorem~3.3, where he proves the generatedness of the anticanonical divisor. Then $|-K_V|$ determines a morphism  
 \[\varphi_{|-K_V|}: V \lra \PN^{g+1},\]
 where $g = \frac{1}{2}(-K_V)^3 +1$ is called the {\em genus} of $V$. Moreover, $\varphi_{|-K_V|}$ is either an embedding, or a 2:1--cover of some smooth variety. By \cite{Isk}, Theorem~7.2., the latter case is very special. The genus is bounded. Iskovskikh shows $2 \le g \le 12$, $g \not= 11$ for $r = 1$ and $g = 5,9,13,17,21$ for $r = 2$. Except for the cases $r = 1$ and $g = 7,9,10,12$ he obtains the description of each Fano threefold as a complete intersection in a (weighted) projective space as listed in the table in the appendix.

\vspace{0.2cm}

Mukai later developed the vector bundle method to classify Fano threefolds. This method leads in particular to a more detailed description in the case of anticanonical embedded Fano threefolds. Our remaining cases $r = 1$ and $g = 7,9,10,12$ are of this type. We have added Mukai's realisation in the table in the appendix for these $4$ cases.  

\vspace{0.2cm}

The reason why we restrict to $\Omega_V^1(1)$ and do not consider higher twists as well, is simply the following. The Euler sequence on projective space $\PN^n$ says 
  \[0 \lra \O_{\PN^n}(-1) \lra \O_{\PN^n}^{\oplus (n+1)} \lra T_{\PN^n}(-1) \lra 0.\]
Applying the functor $\bigwedge^{n-1}$, using $\bigwedge^{n-1}T_{\PN^n}(-1) \simeq \Omega_{\PN^n}^1(2)$, we get a surjection from a sum of $\O_{\PN^n}$'s to $\Omega_{\PN^n}^1(2)$. In this way we see that $\Omega_{\PN^n}^1(2)$ is spanned. If now, for example, $V$ is Fano as above and if $\O_V(1)$ is very ample, then the induced embedding $V \hookrightarrow \PN^n$ defines a map
\[\Omega_{\PN^n}^1(2) \lra \Omega_V^1(2) \lra 0,\]
which shows that $\Omega_V^1(2)$ is spanned as well. By Iskovskikh's classification, $\O_V(1)$ is very ample, except for the cases no.~3,4,8 and 10.  

\begin{proposition} \label{class}
 Let $V$ be an index $r$ Fano threefold of Picard number one and genus $g$. Denote by $\O_V(1) \in \Pic(V)$ an ample generator. The existence of a holomorphic section of $\Omega_V^1(1)$ implies $r = 1$ and $g = 10$ or $12$.
\end{proposition}

We start by proving some general lemmas on the cohomology of twisted $1$--forms, which will later cover all threefolds from the classification.

\begin{lemma} \label{ci}
Let $M$ be a projective manifold of dimension $n \ge 4$. Let $\O_M(1)$ be an ample divisor on $M$ and $V \in |\O_M(d)|$ be a smooth hypersurface. Define $\O_V(1) = \O_M(1)|_V$. Assume $H^0(M, \Omega_M^1(1)) = 0$. If one of the following conditions holds 
\begin{enumerate}
 \item $d \ge 2$,
 \item $d=1$ and $b_2(M) = 1$, 
\end{enumerate}
then $H^0(V, \Omega_V^1(1)) = 0$.
\end{lemma}

\begin{proof}
Write $\O_V(k) = \O_M(k) \otimes \O_V$. The claim follows from standard vanishing theorems applied to the dualized tangent sequence
 \begin{equation} \label{tanseq}
   0 \lra N^*_{V/M} = \O_V(-d) \lra \Omega_M^1|_V \lra \Omega_V^1 \lra 0,
 \end{equation}
and the ideal sequence of $V$ in $M$, tensorized with $\Omega_M^1(1)$:
 \begin{equation} \label{idseq}
   0 \lra \Omega_M^1(1-d) \lra  \Omega_M^1(1) \lra \Omega_M^1(1)|_V \lra 0.
 \end{equation}

1.) Assume $d \ge 2$. Kodaira's vanishing theorem yields $H^i(V, \O_V(1-d)) = 0$ for $i= 0, 1$, so $H^0(V, \Omega_V^1(1)) \simeq H^0(V, \Omega_M^1(1)|_V)$ in (\ref{tanseq}). By Serre duality, $H^1(M, \Omega_M^1(1-d)) \simeq H^{n-1}(M, \Omega_M^{n-1}(d-1))^*$. The latter vanishes by the Kodaira--Akizuki--Nakano vanishing theorem, since $d \ge 2$, $n \ge 3$ and since $\O_M(1)$ is ample. Hence $H^0(M, \Omega_M^1(1)) \lra H^0(V, \Omega_M^1(1)|_V)$ in (\ref{idseq}) is surjective. Since $H^0(M, \Omega_M^1(1)) = 0$ by assumption, we infer $H^0(V, \Omega_M^1(1)|_V) = 0$, and therefore $H^0(V, \Omega_V^1(1)) = 0$.

2.) Assume $d = b_2(M) = 1$. Then $H^0(M, \Omega_M^1) = 0$ in (\ref{idseq}), since by assumption $H^0(M, \Omega_M^1(1)) = 0$. By Lefschetz, $h^0(M, \Omega_M^1) = h^1(M, \O_M) = 0$ implies $h^1(V, \O_V) = 0$. From (\ref{tanseq}), twisted by $\O_V(1)$, we infer 
 \[h^0(V, \Omega_V^1(1)) = h^0(V, \Omega_M^1(1)|_V) - 1 \le h^1(M, \Omega_M^1) - 1.\]
But $h^{1,1}(M) = 1$, since $b_2(M) = 1$, yielding $h^0(V, \Omega_V^1(1)) = 0$.
\end{proof}

The next lemma requires some basic knowledge on weighted projective spaces $\PN(Q) = \PN(q_0, \dots, q_n)$, the {\bf Proj} of $\KC[x_0, \dots, x_n]$, giving $x_i$ weight $q_i$. For details, in particular concerning the definition of the sheaves $\O_{\PN(Q)}(d)$ or $\Omega_{\PN(Q)}^1$, we refer the reader to \cite{Dolg}. Recall that $\PN(Q)$ is called {\em well--formed}, if the $q_i$'s are pairwise relatively prime, and the greatest common divisor of $q_0, \dots, \hat{q_i}, \dots, q_n$ is $1$ for all $i$.

\begin{lemma} \label{weightedps}
 Let $\PN(Q) = \PN(q_0, \dots, q_n)$ be a well--formed weighted projective space for some $n \ge 4$. Let $V \in |\O_{\PN(Q)}(d)|$ be a smooth hypersurface contained in $\PN(Q)_{reg}$, where $d$ is divisible by all the $q_i$'s. Define $\O_V(1) = \O_{\PN(Q)}(1)|_V$. Then $H^0(V, \Omega_V^1(1)) = 0$.
\end{lemma}

\begin{proof}
We first conclude $H^0(\PN(Q), \Omega_{\PN(Q)}^1(1)) = 0$. This follows from the exact Euler sequence on $\PN(Q)$, reading for weighted projective spaces (\cite{Dolg}, \S~2)
 \begin{equation} \label{euler}
  0 \lra \Omega_{\PN(Q)}^1(1) \lra \oplus_{i=0}^n \O_{\PN(Q)}(1-q_i) \stackrel{\rho}{\lra} \O_{\PN(Q)}(1) \lra 0.
 \end{equation}
By \cite{Dolg}, 2.3.4. Corollary, we have $H^j(\PN(Q), \Omega_{\PN(Q)}^l(k)) \not= 0$ only when $j = 0$ and $k > \min_{0 \le i_1 < \dots < i_l \le n}(q_{i_1} + \dots + q_{i_l})$. Hence $H^0(\PN(Q), \Omega_{\PN(Q)}^1(1)) = 0$. 

The sheaf $\O_V(1)$ is free and ample on $V$. We may assume $d \ge 2$, since $d = 1$ implies $q_i = 1$ for all $i$, so $\PN(Q) = \PN^n$, in which case the proof is analogous to the case 1.) of \ref{ci}~Lemma. For $d \ge 2$, since $V$ is supposed to be contained in the smooth locus of $\PN(Q)$, the proof is analogous to 2.) of \ref{ci}~Lemma.
\end{proof}

\begin{proof}[Proof of \ref{class}~Proposition.]
We prove the claim using the classification, for the notation see the table in the appendix. Since $\dim V = 3$, we have $1 \le r \le 4$. By Kobayashi and Ochiai's criterion, if $r = 4$, then $V \simeq \PN^3$, and if $r = 3$, then $V \simeq Q_3$, the quadric hypersurface in $\PN^4$. By Bott's formula, $H^0(\PN^3, \Omega_{\PN^3}^1(1)) = 0$. In the case of the quadric, \ref{ci}~Lemma applies, showing $H^0(Q_3, \Omega_{Q_3}^1(1)) = 0$. It remains to consider the cases $r = 1, 2$.

In the case $r = 2$ we have the following $5$ possibilities: i) $V \in |\O_{\PN(Q)|(6)}$, where $\PN(Q) = \PN(1,1,1,2,3)$ is a weighted projective space; ii) $V \in |\O_{\PN(Q)}(4)|$, $\PN(Q) = \PN(1,1,1,1,2)$; iii) $V$ is a cubic in $\PN^4$; iv) $V \subset \PN^6$ is a complete intersection of two quadrics; v) $V$ is the complete intersection of the $6$--dimensional Grassmannian $Gr(2,5)$ and $3$ hyperplanes in $\PN^9$. In the first two cases i) and ii), \ref{weightedps}~Lemma shows $H^0(V, \Omega_V^1(1)) = 0$. In the latter cases iii) to v) the same is proved by \ref{ci}~Lemma. For v) note, that $H^0(G, \Omega_G^1(1)) = 0$ for a Grassmannian $G$ by \cite{Snow86}, 3.4.~Proposition. Hence, if $r = 2$, then $H^0(V, \Omega_V^1(1)) = 0$.

In the remaining case $r=1$ we have $2 \le g \le 12$, $g \not= 11$ for the genus $g$ of $V$, and we want to prove $g = 10$ or $g = 12$. 

If $g = 2$, then $V \in |\O_{\PN(Q)}(6)|$, $\PN(Q) = \PN(1,1,1,1,3)$, and \ref{weightedps}~Lemma applies showing $H^0(V, \Omega_V^1(1)) = 0$.

If $g = 3$, then $V$ is either a quartic in $\PN^4$, or the following intersection: let $V' \in |\O_{\PN(Q)}(8)|$, where $\PN(Q) = \PN(1,1,1,1,1,4)$, be a general, hence smooth hypersurface. Then $V \in |\O_{V'}(2)|$ is a quartic, where $\O_{V'}(1) = \O_{\PN(Q)}(1)|_{V'}$ by definition. In the first case \ref{ci}~Lemma applies; for the second case apply \ref{weightedps} first, then \ref{ci} for $V \subset V'$ to prove the vanishing $H^0(V, \Omega_V^1(1)) = 0$.

If $g = 4,5$, then $V$ is a complete intersection in some projective space, and if $g = 6,8$, then $V$ is a complete intersection in some Grassmannian. Both cases are clear by \ref{ci}~Lemma and Snow's result on Grassmannians cited above.

If $g = 7,9$, then $V$ is a linear section in the Hermitian symmetric space $M = G/P$ of type DIII for $g = 7$ and CI for $g = 9$ by a result of Mukai (see \cite{Mukai}, \S~2 or \cite{AG5}, \S~5.2.). For the space DIII, $G = SO(10, \KC)$ and $G = Sp(6, \KC)$ in the case CI. The subgroup $P$ of $G$ is maximal parabolic. The cohomolgy of twisted holomorphic forms on manifolds of these types have been studied by Snow in \cite{Snow88}, which gives $H^0(M, \Omega_M^1(1)) = 0$ (see 3.3.~Propsosition and 2.3.~Proposition). 

The only remaining cases are $g = 10$ and $g = 12$ and we are done. 
\end{proof}

\vspace{0.1cm}


\section{Fano threefolds of type $V_{22}$}
\setcounter{equation}{0}

Throughout this section, by $V$ we denote a Fano threefold with Picard number one of genus $12$, i.e. of type $V_{22}$. Then we have natural isomorphisms
 \begin{equation} \label{iso}
    \mbox{$\bigwedge\nolimits^2$}T_V \simeq \Omega_V^1(1) \quad \mbox{ and } \quad \mbox{$\bigwedge^3$}T_V \simeq \O_V(1)
 \end{equation}
and we will sometimes identify these bundles. A general member of the family has a finite automorphism group, hence no vector fields. By \cite{Prokh2}, there are three special types with non--trivial automorphism group: two isolated members $V_{22}^m$ and $V_{22}^a$ with one and two dimensional automorphism group, respectively, and the Mukai--Umemura threefold $V_{22}^s$ with automorphism group $\Sl_2(\KC)$ moving in a one dimensional family. We first show that there are indeed sections in $\Omega_V^1(1)$.

\begin{lemma} \label{ge3}
For $V$ as above of type $V_{22}$, we have $h^0(V, \Omega_V^1(1)) \ge 3$.
\end{lemma}

Now let $\theta \in H^0(V, \Omega_V^1(1))$ be a non--zero section. We may consider $\theta$ as a map $\theta: T_V \to \O_V(1)$. Then $\im(\theta) = \O_V(1) \otimes \I_{\Z(\theta)}$, where $\Z(\theta) \subset V$ is the zero locus of $\theta$. Defining $\F_{\theta} = \ker(\theta)$ we get an exact sequence
 \begin{equation} \label{thetaseq}
  0 \lra \F_{\theta} \lra T_V \stackrel{\theta}{\lra} \O_V(1) \otimes \I_{\Z(\theta)} \lra 0.
 \end{equation}
Since $\O_V(1) \otimes \I_{\Z(\theta)}$ and $\F_{\theta}$ are torsion free, $\F_{\theta}$ is even reflexive (see \cite{OSS}, 1.1.16~Lemma). The generic rank of $\F_{\theta}$ is $2$. Since $\theta$ cannot vanish on a divisor, $\codim(\Z(\theta), V) \ge 2$. Hence $c_1(\F_{\theta}) = 0$.

\vspace{0.2cm}

\begin{proof}[Proof of \ref{ge3}~Lemma.]
A general member $S \in |\O_V(1)|$ is a smooth K3 surface by \cite{Shokurov}. Define $\O_S(1) = \O_V(1)|_S$. The Kodaira--Akizuki--Nakano vanishing theorem implies $h^2(S, \Omega_S^1(1)) = 0$. We will show $H^1(S, \Omega_S^1(1))$ is non--empty: assume to the contrary $h^1(S, \Omega_S^1(1)) = 0$. Then $h^1(S, T_S \otimes N^*_{S/V}) = 0$ by Serre duality, meaning the tangent sequence of $S$ in $V$ splits. This implies
 \[T_V|_S \simeq T_S \oplus \O_S(1).\]
From the ideal sequence we compute $h^0(S, \O_S(1)) = h^0(V, \O_V(1)) -1 = 13$. On the other hand, $h^1(V, T_V(-1)) = h^2(V, \Omega_V^1) = 0$ (see \cite{AG5}, \S~12.2) implies 
 \[h^0(S, T_V|_S) = h^0(V, T_V) - h^0(V, T_V(-1)) \le 3,\] 
since $V$ admits at most $3$ vector fields, a contradiction. Hence $h^1(S, \Omega_S^1(1)) \ge 1$. By Riemann--Roch on $S$, $\chi(S, \Omega_S^1(1)) = \O_V(1)^3 -20 = 2$. We obtain 
 \[h^0(S, \Omega_S^1(1)) = \chi(S, \Omega_S^1(1)) + h^1(S, \Omega_S^1(1)) \ge 3.\] 
From the twisted tangent sequence of $S$ in $V$
 \[0 \lra \O_S \lra \Omega_V^1(1)|_S \lra \Omega_S^1(1) \lra 0\]
we obtain $h^0(S, \Omega_V^1(1)|_S) = h^0(S, \Omega_S^1(1)) + 1 \ge 4$; the sequence
 \[0 \lra \Omega_V^1 \lra \Omega_V^1(1) \lra \Omega_V^1(1)|_S \lra 0\]
then gives $h^0(V, \Omega_V^1(1)) \ge h^0(V, \Omega_V^1(1)|_S) -1 \ge 3$, since $H^1(V, \Omega_V^1) \simeq \KC$.
\end{proof} 

\begin{abs} \label{MUconstr}
{\bf The Mukai--Umemura threefold $V_{22}^s$.} A very special member of the $V_{22}$ family is the almost homogeneous Mukai--Umemura threefold $V_{22}^s$. The construction is as follows (see \cite{MU} for details). Let $M_{12} = \KC[t_0, t_1]_{12}$ be the $\KC$--vector space of homogeneous polynomials in the two variables $t_0, t_1$ of degree $12$. View $M_{12} \simeq \KC^{13}$ as the affine part of $\PN(M_{12} \oplus \KC) \simeq \PN^{13}$ and identify $\PN(M_{12})$ with the hyperplane at infinity. The natural action of $\Sl_2(\KC)$ on $\KC[t_0, t_1]$ induces an action on $\PN(M_{12} \oplus \KC)$. Define
 \[x := t_0t_1(t_0^{10}-11t_0^5t_1^5-t_1^{10}) \in M_{12}.\]
Following Mukai and Umemura, define
 \[V_{22}^s = \overline{\Sl_2(\KC) \cdot [x+1]}.\]
It is not difficult to see that $V_{22}^s$ is indeed a smooth Fano threefold of genus $12$. The action of $\Sl_2(\KC)$ on $V_{22}^s$ has the $3$--dimensional open orbit $O_3 = \Sl_2(\KC)\cdot [x+1]$ and the orbits
 \[O_2 = \Sl_2(\KC) \cdot [t_0t_1^{11}], \quad O_1 = \Sl_2(\KC) \cdot [t_1^{12}]\]
of dimensions $2$ and $1$, respectively. We have $O_1, O_2 \subset \PN(M_{12})$, the hyperplane at infinity. In fact $V_{22}^s = O_1 \cup O_2 \cup O_3$ and $V_{22}^s \cap \PN(M_{12}) = O_1 \cup O_2$, i.e. $O_1 \cup O_2 \in |\O_{V_{22}^s}(1)|$. The orbit $O_2$ is neither open nor closed, $O_1 \simeq \PN^1$ and $\overline{O}_2 = O_1 \cup O_2$. The hyperplane $O_1 \cup O_2 \in |\O_{V_{22}^s}(1)|$ is the hyperplane cut out by lines (cf. \cite{MU}, Lemma~6.1.); it is singular along $O_1 \simeq \PN^1$, the normalization being $\PN^1 \times \PN^1$. This can be seen as follows. Taking a general matrix
 \[\gamma = \left(\begin{array}{cc} a & b \\c & d \end{array}\right) \in \Sl_2(\KC)\]
to compute $O_2$, we find
 \[O_2 = \{[(at_0 + bt_1)(ct_0 + dt_1)^{11}] \mid ad-bc = 1\} \subset \PN(M_{12}).\]
The map $\nu: \PN^1 \times \PN^1 \to \PN(M_{12})$ defined by $[a:b] \times [c:d] \mapsto [(at_0 + bt_1)(ct_0 + dt_1)^{11}]$, i.e. by a subsystem of $|\O_{\PN^1\times\PN^1}(1,11)|$, is then a normalization map of $\overline{O}_2$. Here $\nu$ is equivariant with respect to the action on $\overline{O}_2$ and the transposed diagonal action on $\PN^1 \times \PN^1$, i.e. $\nu \gamma^t = \gamma \nu$ for any $\gamma \in \Sl_2(\KC)$. The nonnormal locus of $\overline{O}_2 = O_1 \cup O_2$ is $\nu(\Delta) = O_1$, where $\Delta$ denotes the diagonal in $\PN^1 \times \PN^1$. We see from this description that $O_1 \cup O_2$ is indeed cut out by lines.

The equivariance of $\nu$ implies the following: we have a map
 \[H^0(\Delta, T_{\Delta}) \stackrel{i}{\hookrightarrow} H^0(\PN^1 \times \PN^1, T_{\PN^1 \times \PN^1})\]
defined as follows. For $X \in H^0(\Delta, T_{\Delta})$ define $i(X)(p,q) = (X(p), X(q)) \in p_1^*T_{\PN^1} \oplus p_2^*T_{\PN^1} = T_{\PN^1 \times \PN^1}$, where $p_i$ denote the projections. Then for any $Y \in H^0(V_{22}^s, T_{V_{22}^s})$ we have
 \[\nu^*Y \in \im \big(H^0(\Delta, T_{\Delta}) \lra H^0(\PN^1 \times \PN^1, T_{\PN^1 \times \PN^1}) \lra H^0(\PN^1 \times \PN^1, \varphi^*T_{V_{22}^s})\big).\]
\end{abs}

\begin{proposition} \label{MU}
\begin{enumerate}
  \item Let $X,Y \in H^0(V_{22}^s, T_{V_{22}^s}) \simeq {\mathfrak sl}_2(\KC)$ be linearly independent vector fields and define $\theta_{X,Y} = X \wedge Y \in H^0(V_{22}^s, \bigwedge^2T_{V_{22}^s})$. Then 
  \begin{enumerate}
    \item $\Z(\theta_{X, Y})_{red} = O_1 \cup (\mbox{rational curve}) \subset O_2$, 
    \item $\Z(d\theta_{X,Y} \wedge \theta_{X,Y}) = V$ or $O_1 \cup O_2$, depending on whether $X$ and $Y$ generate a subalgebra of ${\mathfrak sl}_2(\KC)$ or not.
    \item $\F_{\theta_{X,Y}} \simeq \O_{V_{22}^s}^{\oplus 2}$, i.e. we have the exact sequence
 \[0 \lra \O_{V_{22}^s}^{\oplus 2} \lra T_{V_{22}^s} \stackrel{\theta_{X,Y}}{\lra} \O_{V_{22}^s}(1) \otimes \I_{\Z(\theta_{X,Y})} \lra 0.\]
  \end{enumerate}
 \item $H^0(V_{22}^s, \bigwedge^2T_{V_{22}^s}) \simeq \bigwedge^2H^0(V_{22}^s, T_{V_{22}^s}) \simeq \KC^3$, meaning that any section in $H^0(V_{22}^s, \Omega_{V_{22}^s}^1(1))$ is as in 1.).
\end{enumerate}
\end{proposition}

\begin{remark}
  In 1.2.), if $X, Y \in H^0(V_{22}^s, T_{V_{22}^s}) \simeq {\mathfrak sl}_2(\KC)$ define a subalgebra of ${\mathfrak sl}_2(\KC)$, then $d\theta_{X, Y} \wedge \theta_{X, Y} \equiv 0$, and we have a foliation. The leaves are the orbits of the corresponding subgroup of $\Sl_2(\KC)$. In general, however, we will have $d\theta_{X, Y} \wedge \theta_{X, Y} \not\equiv 0$, and $\Z(d\theta_{X,Y} \wedge \theta_{X,Y}) = O_1 \cup O_2$.
\end{remark}

\begin{proof}[Proof of \ref{MU}~Proposition.]
We write $V$ instead of $V_{22}^s$ for simplicity.

1.1.) Let $X,Y \in H^0(V, T_V)$ be two linearly independent vector fields. Note $H^0(V, T_V) = \KC^3$. Using (\ref{iso}), we may think of $X \wedge Y$ as a section of $\Omega_V^1(1)$. The zero set of this section is $\Z = \{p \in V \mid (X \wedge Y)(p) = 0\}$. We know $\dim_{\KC} \Z \le 1$. Since $T_V|_{O_3}$ is generated by three sections, $O_3 \cap \Z = \emptyset$. Hence, set theoretically, $\Z \subset O_1 \cup O_2$. From above:
 \[\nu^*X, \nu^*Y \in \im \big(H^0(\Delta, T_{\Delta}) \lra H^0(\PN^1 \times \PN^1, \nu^*T_V)\big).\]
It is then clear from this description that $\Delta$ is part of the zero locus of $\nu^*(X \wedge Y)$. It is moreover clear that
 \[\nu^*(X\wedge Y) \in \im \big(H^0(\PN^1 \times \PN^1, \mbox{$\bigwedge^2$}T_{\PN^1 \times \PN^1}) \lra H^0(\PN^1 \times \PN^1, \nu^*\mbox{$\bigwedge^2$}T_V)\big).\]
From $\bigwedge^2T_{\PN^1 \times \PN^1} = p_1^*\O_{\PN^1}(2) \otimes p_2^*\O_{\PN^1}(2)$ we infer the vanishing locus of $\nu^*(X \wedge Y)$ is either $2\Delta$ or $\Delta + \Delta'$, where $\Delta' \in |\O_{\PN^1 \times \PN^1}(1,1)|$. In the first case, set theoretically, $\Z = \nu(2\Delta)$, in the latter case $\Z = \nu(\Delta) \cup \nu(\mbox{rational curve of degree 12})$. In any case, $\Z_{red} = O_1 \cup (\mbox{rational curve})$. This proves 1.1.).

\vspace{0.2cm}

1.3.) and 2.). Define $W = \{X \wedge Y \mid X,Y \in H^0(V, T_V)\} \subset H^0(V, \bigwedge^2T_V)$. Three generating vector fields in $H^0(V, T_V)$ are pairwise independent on $O_3$, implying $W \simeq \bigwedge^2H^0(V, T_V)$, a three dimensional vector space. We want to show $W = H^0(V, \Omega_V^1(1))$. Using the notation from (\ref{thetaseq}) we prove the equivalences
  \[\theta \in W\backslash\{0\} \quad \Longleftrightarrow \quad h^0(V, \F_{\theta}) \ge 2 \quad \Longleftrightarrow \quad \F_{\theta} \simeq \O_V^{\oplus 2}.\]
The equivalences imply 1.3.).

We first prove $h^0(V, \F_{\theta}) \ge 2$ implies $\F_{\theta} \simeq \O_V^{\oplus 2}$. From $\F_{\theta} \hookrightarrow T_V$ we infer $h^0(V, \F_{\theta}) \le 3$. Three vector fields generate $T_V$ on $O_3$. Then they cannot be all contained in $H^0(V, \F_{\theta})$, since $\F_{\theta}$ is generically of rank two. Hence $h^0(V, \F_{\theta}) = 2$. Let $X_0, Y_0 \in H^0(V, T_V)$ be generators of $H^0(V, \F_{\theta})$, i.e. $\theta(X_0) = \theta(Y_0) = 0$. Define $\Z_0 = \{p \in V \mid (X_0 \wedge Y_0)(p) = 0\}$. Then $\codim(\Z_0, V) = 2$, since $X_0 \wedge Y_0$ vanishes on a curve by 1.). This gives a map $\O_V^{\oplus 2} \to \F_{\theta}$, which is surjective away from $\Z_0$. This shows $\O_V^{\oplus 2} \simeq \F_{\theta}$, since $\F_{\theta}$ is reflexive and $c_1(\F_{\theta}) = 0$. 

Now assume $\F_{\theta} \simeq \O_V^{\oplus 2}$. We prove that then $\theta \in W\backslash\{0\}$. Indeed, using the notation from above, we may assume $H^0(V, \F_{\theta})$ is generated by two vector fields $X_0,Y_0$. By construction, the map $i:\F_{\theta} \simeq \O_V^{\oplus 2} \hookrightarrow T_V$ is then defined by $(f,g) \mapsto fX_0 + gY_0$. Consider on the other hand $\theta_0: T_V \to \O_V(1) \otimes \I_{\Z_0}$ defined by $X_0 \wedge Y_0$. Denote the kernel by $\F_0$. Then $\F_0 \simeq \O_V^{\oplus 2}$ as above, and the inclusion $\F_0 \simeq \O_V^{\oplus 2} \hookrightarrow T_V$ is the same map as $i$. Therefore the cokernel maps must coincide, meaning $\theta = \lambda X_0 \wedge Y_0$ for some $\lambda \in \KC^*$ (and $\Z(\theta) = \Z_0$).

Finally assume $0 \not\equiv \theta \in W$. Then $\F_{\theta} \simeq \O_V^{\oplus 2}$ as above, hence $h^0(V, \F_{\theta}) = 2$.

\vspace{0.2cm}

To finally prove $W = H^0(V, \Omega_V^1(1))$, consider some $\theta_0 \in W\backslash\{0\}$. Then $\F_{\theta_0} \simeq \O_V^{\oplus 2}$, as we have seen. Since $\O_V$ is rigid, for $\theta_t$ chosen from some (analytically) open neighborhood $U(\theta_0) \subset H^0(V, \Omega_V^1(1))$ of $\theta_0$, we also have $\F_{\theta_t} \simeq \O_V^{\oplus 2}$. The above equivalences show $U(\theta_0) \subset W$, implying $W = H^0(V, \Omega_V^1(1))$. Point 2.) is proved.

\vspace{0.2cm}

1.2.) To determine $d\theta_{X,Y} \wedge \theta_{X,Y}$, consider the map
 \[\O_V \simeq \mbox{$\bigwedge^2$}\F_{X,Y} \lra \O_V(1)\]
induced by $\theta_{X,Y} \circ [-,-] = X \wedge Y \wedge [-,-]$. We see that the zero set of $d\theta_{X,Y} \wedge \theta_{X,Y}$ is the zero set of $X \wedge Y \wedge [X,Y]$, with (\ref{iso}) viewed as a section of $\O_V(1)$. If $[X,Y]\in \langle X, Y\rangle_{\KC}$, then $\Z(d\theta_{X,Y} \wedge \theta_{X,Y}) = V$. Otherwise, choose $Z$ such that $H^0(V, T_V) = \langle X, Y, Z\rangle_{\KC}$. We have to find the zero set of $X \wedge Y \wedge Z$. On $O_3$, the three sections are independent, so they define a nonzero section of $\O_V(1)$, vanishing on the complement of $O_3$. We finally conclude $\Z(d\theta_{X,Y} \wedge \theta_{X,Y}) = O_1 \cup O_2 \in |\O_V(1)|$.
\end{proof}

\begin{abs} \label{V22family}
{\bf Family of Fano threefolds of type $V_{22}$.} By Mukai's construction (see \cite{Mukai}, or \cite{AG5}, \S 5.2.), any Fano threefold $V$ of type $V_{22}$ can be embedded into the Grassmannian ${\rm Gr}(7,3)$ of $3$--dimensional quotient spaces of $\KC^7$. Let ${\mathcal Q}$ be the universal quotient bundle on the Grassmannian. Then $V$ is defined as zero locus of $3$ sections in $\bigwedge^2{\mathcal Q}$. The parameter space of $V_{22}$ is birationally equivalent to the moduli space of curves of genus $3$ by \cite{EPS} or \cite{AG5}, p.114, hence $6$--dimensional and irreducible. Assume that $V$ is not the Mukai--Umemura threefold. Then the divisor cut out by lines is a reduced, irreducible divisor from $|\O_V(2)|$ (see \cite{AG5}, \S 4.2, \cite{Prokh} and \cite{SchIl}), and the splitting type of $T_V$ on a general line is $(2,0,-1)$. 
\end{abs}

\begin{proposition} \label{allgV22}
Let $V$ be general of type $V_{22}$. Then $h^0(V, \Omega_V^1(1)) = 3$ and $d\theta \wedge \theta \equiv 0$ for any $\theta \in H^0(V, \Omega_V^1(1))$.
\end{proposition}

\begin{proof}
We will apply semicontinuity on the family of Fano threefolds of type $V_{22}$. Let $V$ be a general member and $V^s = V_{22}^s$ the Mukai--Umemura threefold, a special member. Then
 \[h^0(V, \Omega_V^1(1)) \le h^0(V^s, \Omega_{V^s}^1(1)) = 3,\]
by \ref{MU}~Proposition. On the other hand $h^0(V, \Omega_V^1(1)) \ge 3$ by \ref{ge3}~Lemma, showing $h^0(V, \Omega_V^1(1)) = 3$.

Let $\theta \in H^0(V, \Omega_V^1(1))$ be a non--zero section. We want to prove $d\theta \wedge \theta \equiv 0$. Since $H^0(V, \Omega_V^1(1))$ is threedimensional, $\theta$ is a deformation of some $\theta_s \in H^0(V^s, \Omega_{V^s}^1(1))$. Define the kernels $\F_{\theta}$ and $\F_{\theta_s}$ as in (\ref{thetaseq}). By \ref{MU}~Proposition, $\F_{\theta_s} \simeq \O_{V^s}^{\oplus 2}$. 

On $V$ we have the exact sequence
 \begin{equation} \label{theta}
   0 \lra \F_{\theta} \lra T_V \lra \O_V(1) \otimes \I_{\Z(\theta)} \lra 0.
 \end{equation}
We will show that $\theta$ vanishes in more than one point on a general line $l \subset V$. First, we may assume that $l$ does not meet the codimension $3$ locus, where $\F_{\theta}$ is not free. Therefore $\F_{\theta}|_l$ is a rank two vector bundle of degree $0$. Let $l_s$ be a line in $V^s$, obtained by deforming $l$. By semicontinuity, $h^0(l, \F_{\theta}(-1)|_l) \le h^0(l_0, \F_{\theta_s}(-1)|_{l_s}) = 0$. This shows $\F_{\theta}|_l \simeq \O_l^{\oplus 2}$. 

The splitting type of $T_V$ on $l$ is $T_V|_l = \O_l(2) \oplus \O_l \oplus \O_l(-1)$ (\cite{AG5}, Theorem~4.2.7). The restriction $\I_{\Z(\theta)} \otimes \O_l$ might not be torsion free, but nevertheless, the vanishing order of $\theta$ on $l$ is exactly the (negative) degree of the free part, since $\Z(\theta)$ meets $l$ only in points. The restriction of (\ref{theta}) hence looks like
 \[0 \lra \O_l^{\oplus 2} \stackrel{\alpha}{\lra}  \O_l(2) \oplus \O_l \oplus \O_l(-1) \lra \O_l(-a+1) \oplus \tau \lra 0,\]
where $\tau$ is a torsion sheaf, and $a$ is the order of $\Z(\theta) \cap l$ we are looking for. Computing $H^1$, we find $a = 2$. We have proved, that $\theta$ vanishes in $2$ points on $l$.

Consider now $d\theta \wedge \theta$. Since $d\theta \wedge \theta$ obviously vanishes in the zeroes of $\theta$, it vanishes in two points on a general line $l$. Since $d\theta \wedge \theta \in |\O_V(1)|$, it follows $d\theta \wedge \theta|_l \equiv 0$. This implies, that $d\theta \wedge \theta$ vanishes on the whole divisor cut out by lines, which is an element in $|\O_V(2)|$, if $V \not= V_{22}^s$. This shows $d\theta \wedge \theta \equiv 0$.
\end{proof}

\vspace{0.1cm}


\section{Fano threefolds of type $V_{18}$} \label{sec v18}
\setcounter{equation}{0}

Let $M$ be the $5$ dimensional contact manifold, homogeneous under the exceptional group $G_2$. Naturally embedded in $\PN^{13}$, the contact bundle of $M$ is the fundamental divisor $L = \O_M(1) = \O_{\PN^{13}}(1)|_M$. We use
\[0 \lra F \lra T_M \stackrel{\theta_M}{\lra} L \lra 0\]
to describe the contact sequence. The contact form $\theta_M \in H^0(M, \Omega_M^1(1))$ is unique up to multiples.

By Mukai's construction, a Fano threefold $V$ of type $V_{18}$ is a complete intersection of two hyperplanes $H_1, H_2 \in |\O_M(1)|$ in our contact manifold $M$. We do not have vector fields on $V$. With this interpretation of $V$, we first prove 

\begin{proposition}
  For $V$ of type $V_{18}$ we have $H^0(V, \Omega_V^1(1)) = \KC$, a generating section being the image of $\theta_M$ under $H^0(M, \Omega_M^1(1)) \lra H^0(V, \Omega_V^1(1))$. 
\end{proposition}

\begin{proof}
Since $-d\theta_M = \theta_M([-,-]): F \times F \lra L$ is non--degenerate, Frobenius theorem implies that if $W$ is a submanifold of $M$ and $T_W \subset F|_W$, then $\dim W < 3$. Then $T_V$ cannot be contained in $F|_V$, and from 
 \[\xymatrix{ && T_V \ar@{^{(}->}[d] \ar@{..>}[dr]^{\theta} &&\\
     0 \ar[r] & F|_V \ar[r] & T_M|_V \ar[r]^{\theta_M|_V} & \O_V(1) \ar[r] & 0,}\]
we see that $\theta_M$ is mapped to a non--vanishing section $\theta$ of $\Omega_V^1(1)$ under the natural map $\Omega_M^1(1) \to \Omega_V^1(1)$. Analogously we see that $\theta_M$ induces a non--vanishing section of $\Omega_{H_1}^1(1)$. 

To show $H^0(V, \Omega_V^1(1)) = \KC$, we use the dualized tangent sequence of $V$ in $H_1$ and the ideal sequence. The first is
 \[0 \lra \O_V \lra \Omega_{H_1}^1(1)|_V \lra \Omega_V^1(1) \lra 0,\]
yielding $h^0(V, \Omega_V^1(1)) = h^0(V, \Omega_{H_1}^1(1)|_V) -1$. By adjunction formula and Lefschetz, $H_1$ is a Fano manifold of Picard number one and $h^1(H_1, \Omega_{H_1}^1) = 1$. The ideal sequence, tensorized with $\Omega_{H_1}^1(1)$ reads
 \[0 \lra \Omega_{H_1}^1 \lra \Omega_{H_1}^1(1) \lra \Omega_{H_1}^1(1)|_V \lra 0,\]
and we get $h^0(V, \Omega_{H_1}^1(1)|_V) \le h^0(H_1, \Omega_{H_1}^1(1)) + 1$. The two estimations yield $h^0(V, \Omega_V^1(1)) \le h^0(H_1, \Omega_{H_1}^1(1))$. 

Analogously, using the same sequences for $H_1$ in $M$, we find $h^0(H_1, \Omega_{H_1}^1(1)) \le h^0(M, \Omega_M^1(1)) = 1$, and we conclude $h^0(V, \Omega_V^1(1)) \le 1$. 
\end{proof}

\vspace{0.2cm}
 
To describe its zero locus as well as $d\theta \wedge \theta \in H^0(V, \O_V(1))$, we now briefly recall the group theoretic background of $M$ and its contact structure. We refer to \cite{Beau} for details.

Instead of considering merely the exceptional group $G_2$, we study an arbitrary simple complex Lie group $G$. Let $\g$ be its Lie algebra. Note that $\g$ and $\g^*$ are isomorphic via the Cartan killing form $\langle-,-\rangle$ (and because of this we will sometimes write $\g$ where perhaps $\g^*$ would be more apropriate in the sequel). There exists exactly one closed orbit $M$ of the adjoint action of $G$ on $\PN(\g)$. Let $L = \O_{\PN(\g)}(1)|_M$. We briefly sketch the idea of the following well known result: {\em $M$ carries a contact structure with contact line bundle $L$ if and only if the dimension of $M$ is odd.}

One direction is trivial. Indeed, if $M$ carries a contact structure $\theta_M$ with contact line bundle $L$, then the pull back of $\theta_M$ to the total space of $L$ induces a symplectic structure on $L$, showing that $\dim M$ must be odd. To prove that the convers holds in the above situation, we first define this symplectic structure, before showing that it comes from a contact from.

Let $M^{\circ}$ be the orbit of $G$ under the adjoined action of $G$ on $\g$, such that $\PN(M^{\circ}) = M$. The tangent space $T_{M^{\circ}}(Z)$ is canonically isomorphic to $\g/\ker\ad(Z)$ for any $Z \in M^{\circ}$. On $M^{\circ}$ we have a nowhere degenerated symplectic form, locally defined by
 \begin{equation} \label{KoKi}
   \omega_{Z}: \quad T_{M^{\circ}}(Z) \times T_{M^{\circ}}(Z) \lra \KC, \quad (X, Y) \mapsto \langle[X, Y], Z\rangle,
 \end{equation}
which is nothing but the Kostant--Kirillov symplectic form, usually rather defined via the coadjoined representation. Note that $\omega_{Z}$ is well defined at $Z$ by Jacobi's formula. The existence of $\omega_{Z}$ implies that $\dim M^{\circ}$ is even. Now assume $\dim M$ is odd. 

The dimension dropping by one, by going from $M^{\circ}$ to $M$, means $M^{\circ}$ is the total space of $L$ over $M$. This is the case if and only if $Z \in M^{\circ}$ implies $cZ \in M^{\circ}$ for any $Z \in M^{\circ}$ and $c \in \KC^*$. But $M^{\circ}$ is an orbit, so this is the case if and only if $Z$ and $cZ$ are conjugated under the adjoined action for any choice of $c \in \KC^*$ and $Z \in M^{\circ}$. This holds if and only if for any $Z \in M^{\circ}$ there exists an $H_Z \in \g$ such that $[H_Z, Z] = Z$.

The existence of an $H_Z \in \g$ for any $Z \in M^{\circ}$ such that $[H_Z, Z] = Z$ implies that for any $Z \in M^{\circ}$ we have $\z_{[Z]} \subset Z^{\perp}$, where 
 \begin{equation} \label{zz}
   \z_{[Z]} = \{X \in \g \mid [X, Z] = \lambda Z \mbox{ for some } \lambda \in \KC\}
 \end{equation}
and $Z^{\perp} = \{X \in \g \mid \langle X, Z\rangle = 0\}$. Indeed, if $X \in  \z_{[Z]}$ and $[X, Z] = \lambda Z$, $\lambda \not= 0$, then $\langle X, Z\rangle = \lambda^{-1} \langle X, [X, Z]\rangle = 0$. If $X \in \z_{[Z]}$ and $[X, Z] = 0$, we pick $H_Z$ from above satisfying $[H_Z, Z] = Z$, and we see $\langle X, Z\rangle = \langle X, [H_Z, Z]\rangle = \langle [X, Z], H_Z\rangle = 0$.

As in the case of $T_{M^{\circ}}(Z)$ we have a canonical isomorphism for the tangent space $T_M([Z])$ of $M$ at $[Z] \in M$, $Z \in \g$:
  \[T_M({[Z]}) \simeq \g/\z_{[Z]}.\]
At the point $[Z] \in M$, the total space of $\O_{\PN(\g)}(1)$ is isomorphic to $\g/Z^{\perp}$ (using $\g \simeq \g^*$). Since $\z_{[Z]} \subset Z^{\perp}$, we have a well defined surjection $T_M({[Z]}) \to L([Z])$, which glues, yielding a bundle sequence
  \[0 \lra F \lra T_M \stackrel{\theta_M}{\lra} L \lra 0.\]
By construction, the pull back of the contact form to the total space of $L$ is the Kostant--Kirillov form (\ref{KoKi}), showing that $\theta_M$ indeed defines a contact structure. Alternatively one may consider the induced map
\[-d\theta_M = \theta_M([-,-]): F \times F \lra L, \quad \mbox{given by} \quad (X, Y) \mapsto [X, Y] \mod Z^{\perp}.\]
This map is non--degenerate. Indeed, at $[Z] \in M$ we have $F = Z^{\perp}/\z_{[Z]}$. Fix some $Y \in Z^{\perp}$ and assume $[X, Y] \in Z^{\perp}$ for any $X \in Z^{\perp}$. Then $\langle [X, Y], Z\rangle = 0$, implying that the hyperplane $\langle -, [Y, Z]\rangle = 0$ contains the hyperplane $Z^{\perp}$. Then $[Y, Z] = \lambda Z$ for some $\lambda \in \KC$. Then $Y \in \z_{[Z]}$, showing that the map is indeed non--degenerate.

\vspace{0.2cm}

The construction of the homogeneous contact manifold $M$ shows that the contact sequence induces an isomorphism $H^0(M, T_M) \simeq H^0(M, L)$. Hyperplane sections of $M$ and vector fields may in this way be identified. Assume from now on that $V$ is cut out by the two smooth general hyperplanes $H_1$ and $H_2$, which are in this sense given by the two vector fields 
  \[X_1, X_2 \in H^0(M, T_M).\]
In our situation $H^0(M, T_M) = \g$, so we may think of $X_1, X_2$ as elements of $\g$. In explicit form, $H_i$ is now given by $\{[Z] \in M \mid \langle Z, X_i\rangle = 0\}$ and
  \[V = \{[Z] \in M \mid \langle Z, X_i\rangle = 0, \mbox{ for }i = 1, 2\}.\]
Since $V$ is again a Fano manifold of index $1$, we have again canonical isomorphisms 
 \[\mbox{$\bigwedge\nolimits^2$}T_V \simeq \Omega_V^1(1) \quad \mbox{ and } \quad \mbox{$\bigwedge^3$}T_V \simeq \O_V(1)\]
With this description, we can interpret $\theta$ as follows:

\begin{proposition}
  Let $V$ be Fano of type $V_{18}$, given as above as a complete intersection of hyperplanes $H_1, H_2$ of the homogeneous $G_2$--contact manifold $M$, induced by vector fields $X_1, X_2 \in \g_2$ on $M$. Then
\[\theta_{X_1, X_2} = X_1 \wedge X_2|_V \in H^0(V, \mbox{$\bigwedge^2$}T_V)\]
is non--vanishing and may be thought of as the pull back of the contact structure $\theta_M$.  The vanishing locus of $d\theta_{X_1, X_2} \wedge \theta_{X_1, X_2}$ is the vanishing locus of $[X_1, X_2]|_V \in H^0(V, \O_V(1))$.
\end{proposition}

\begin{proof}
We begin with considering a single smooth general hyperplane section $H_1$ of $M$, cut out by a section corresponding to $X_1 \in H^0(M, T_M) = \g_2$, i.e.,
  \[H_1 = \{[Z] \in M \mid \langle Z, X_1\rangle = 0\}\]
as above. The contact form induces a nonzero section $\theta_{H_1} \in H^0(H_1, \Omega_{H_1}^1(1))$. We are interested in finding the points where $\theta_{H_1}: T_{H_1} \to \O_{H_1}(1)$ drops rank. The tangent space of the hyperplane $H_1$ at a point $[Z] \in H_1$ has the following canonical description 
\begin{equation} \label{TH}
   T_{H_1}([Z]) = [X_1, Z]^{\perp}/\z_{[Z]}
 \end{equation} 
with $\z_{[Z]}$ as in (\ref{zz}). Note that $\z_{[Z]} \subset [X_1, Z]^{\perp}$ for $[Z] \in H_1$. From this description we see: $X_1$ viewed as a vector field on $M$ is contained in $T_{H_1}([Z])$ for every $[Z] \in H_1$, implying 
  \[X_1 \in H^0(H_1, T_{H_1}).\]

The form $\theta_{H_1}$ drops rank preciseley at those points $[Z] \in H_1$, where the contact bundle $F$ and $T_{H_1}$ define the same hyperplane of $T_M$. Hence $\theta_{H_1}$ drops rank precisely at those $[Z] \in H_1$ where $[X_1,Z]^{\perp} = Z^{\perp}$, which in turn holds precisely for those $[Z] \in M$ satisfying $[X_1, Z] = \lambda Z$ for some $\lambda \in \KC^*$. The latter condition is equivalent to $X_1 \in \z_{[Z]}$. For the equivalence note again that $\z_{[Z]} \subset Z^{\perp}$ and that $H_1$ is smooth.

If we view $X_1$ as a vector field of $M$, then those points $[Z] \in M$, where $X_1 \in \z_{[Z]}$, are the zeroes of $X_1$. We have proved
  \[\mbox{Zero locus of } \theta_{H_1} \in H^0(H_1, \Omega_{H_1}^1(1)) = \mbox{Zero locus of } X_1 \in H^0(H_1, T_{H_1}).\]
The tangent bundle of homogeneous $M$ is globally generated. The vanishing locus of a general section consists of points, their number being equal to the highest Chern class of $M$. We finally conclude (see \ref{c5M}~Lemma):
  \[\mbox{Zero locus of } \theta_{H_1} \in H^0(H_1, \Omega_{H_1}^1(1)) = c_4(T_{H_1}) = c_5(T_M)= 6 \mbox{ points}.\]
Purely in terms of Chern classes, our result contains (and shows) the following equality of Chern classes, which also follow from the tangent sequence combined with the contact sequence: $c_4(T_{H_1}) = c_4(F|_{H_1}) = c_4(F^*(1)|_{H_1}) = c_4(\Omega_{H_1}^1(1))$ and $c_5(T_M) = c_4(T_{H_1})$.

Concerning $d\theta_{H_1} \wedge \theta_{H_1} \in H^0(H_1, \bigwedge^3 \Omega_{H_1}^1 \otimes \O_{H_1}(2))$. Since $\bigwedge^3 \Omega_{H_1}^1 \otimes \O_{H_1}(2) \simeq T_{H_1}$, we may view $d\theta_{H_1} \wedge \theta_{H_1}$ as a vector field on $H_1$. Writing down an explicit isomorphism, we find that the vanishing locus of $d\theta_{H_1} \wedge \theta_{H_1}$ coincides with the vanishing locus of $X_1 \in H^0(H_1, T_{H_1})$, which are $6$ points.

\vspace{0.2cm}

Now consider $V$ from above, the complete intersection of the hyperplanes $H_1, H_2$ corresponding to $X_1, X_2 \in H^0(M, T_M)$. First, since $X_i$ is in the kernel of the map $T_M|_{H_i} \lra \O(1)$ on global sections, it is clear from
\[0 \lra \mbox{$\bigwedge^2$}T_V \lra \mbox{$\bigwedge^2$}T_{H_i}|_V \lra T_V(1) \lra 0\]  
\[0 \lra \mbox{$\bigwedge^2$}T_{H_i} \lra \mbox{$\bigwedge^2$}T_M|_{H_i} \lra T_{H_i}(1) \lra 0\]
that $X_1 \wedge X_2$ indeed defines a nonzero section of $\bigwedge^2T_V$. Since $H^0(V, \Omega_V^1(1))$ is one dimensional, we may take this section to be the pull back of $\theta_M$. Alternatively, one may derive a description of $T_V([Z])$ and conclude as above, that the pull back of $\theta_M$ on $V$ drops rank precisely at the vanishing points of $X_1 \wedge X_2$.

We use the following identifications do determine $d\theta \wedge \theta$:
\[\xymatrix{0 \ar[r] & \mbox{$\bigwedge^2$}\Omega_V^1 \otimes \O_V(1) \ar[r] \ar@{=}[d] & \mbox{$\bigwedge^3$}\Omega_{H_1}^1 \otimes \O_{H_1}(2)|_V \ar[r] \ar@{=}[d] & \mbox{$\bigwedge^3$}\Omega_V^1 \otimes \O_V(2) \ar[r] \ar@{=}[d] & 0\\
0 \ar[r] & T_V \ar[r] & T_{H_1}|_V \ar[r] & \O_V(1) \ar[r] & 0}\]
The pull back of $d\theta_M \wedge \theta_M$ to $H_1$, using the identification $\bigwedge^3\Omega_{H_1}^1 \otimes \O_{H_1}(2) = T_{H_1}$, yields $X_1$ as we saw above. The image of $X_1$ under $H^0(H_1, T_{H_1}) \lra H^0(V, \O_V(1))$ is the section induced by the vector field $[X_1, X_2]$. This is clear from the identification $\bigwedge^2 \Omega_V^1 \otimes \O_V(1) = T_V$ and, for example, the pointwise description of the tangent bundle on $V$, analogous to (\ref{TH}).

Since $X_1$ and $X_2$ were chosen general, it is clear, that $[X_1,X_2]$ does not vanish, and is different from $X_1$ and $X_2$  as elements in $H^0(M, T_M) = \g_2$. It therefore defines a non--zero section in $\O_V(1) = \O_M(1)|_V$.

\vspace{0.2cm}

We check for entertainment that the zero locus of $X_1 \wedge X_2$, viewed as a section of $\bigwedge^2T_V$, is indeed contained in the vanishing locus of $[X_1, X_2]$, viewed as a section of $\O_V(1)$. This must necessarily be the case, since $\Z(\theta) \subset \Z(d\theta \wedge \theta)$. 

If the wedge product $X_1 \wedge X_2$ vanishes at a point $[Z] \in V$, then $X_1$ and $X_2$, evaluated at $[Z]$, are dependent, meaning $\lambda Z = \lambda_1[X_1, Z] + \lambda_2[X_2, Z]$ for some $\lambda, \lambda_1, \lambda_2 \in \KC$. We may assume $\lambda_1 \not= 0$. Applying $\langle -, X_2\rangle$, we find on the left hand side $\lambda \langle Z, X_2\rangle$. Since $[Z]$ is a point on $V$, this is zero. The right hand side then reads $\lambda_1\langle [X_1, Z], X_2\rangle + \lambda_2\langle [X_2, Z], X_2\rangle = \lambda_1\langle [X_1, Z], X_2\rangle= 0$, using $[X_2, X_2] = 0$. Since $\lambda_1\not= 0$, we conclude $\langle [X_1, X_2], Z\rangle = 0$, as desired.
\end{proof}


\begin{lemma} \label{c5M}
 Let $M$ be a quotient of the exceptional simple Lie group $G_2$ by a maximal parabolic subgroup. Then $\dim M = 5$ and $c_5(T_M) = 6$.
\end{lemma}

\begin{proof}
Let $G$ be the exceptional group of type $G_2$ and $B \subset G$ a Borel group. Then there are two (maximal) parabolic subgroups $P_1, P_2$ in $G$ containing $B$. The corresponding homogeneous manifolds are $M_1 = G/P_1$, a $5$--dimensional quadric, and $M_2 = G/P_2$, the $5$--dimensional contact manifold associated to $G$. The dimension of the homogeneous manifold $M_B = G/B$ is $6$. We will show that the highest Chern classes of $M_1$ and $M_2$ coincide. We have the following diagram:
 \[\xymatrix{M_B \ar[r]^{\pi_1} \ar[d]^{\pi_2} & M_1\\
     M_2 &}\]
Let ${\LB}_1$, $\LB_2$ be the (globally generated) fundamental line bundles on $M_B$. Then $-K_{M_B} = 2{\LB}_1 + 2{\LB}_2$ and ${\LB}_i = \pi_i^*L_i$, where $L_i$ is the fundamental line bundle on $M_i$. We have $-K_{M_1} = 5L_1$ and $-K_{M_2} = 3L_2$. All these facts can be found for example in \cite{Akh}. The fibers $F_i = P_i/B$ of $\pi_i$ are so--called $\alpha$--lines in $M_B$, that are smooth rational curves with the property ${\LB}_i.F_j = 1$ for $i \not= j$. This can be easily checked or can be found for example in \cite{Kollar}. In particular, the projections $\pi_i$ are $\PN^1$--bundles. Consider the relative tangent sequences
 \[0 \lra T_{M_B/M_i} \lra T_{M_B} \lra \pi^*_iT_{M_i} \lra 0\]
for $i = 1,2$. The realtive tangent bundles $T_{M_B/M_i}$ are line bundles, which we get by computing the determinant of the above sequence: $T_{M_B/M_1} = -3{\LB}_1 + 2{\LB}_2$ and $T_{M_B/M_2} = 2{\LB}_1 -{\LB}_2$. Since Chern polynomials in short exact sequences are multiplicative, we have
 \[c_t(\pi_1^*T_{M_1}).(1 + c_1(T_{M_B/M_1})t) = c_t(\pi_2^*T_{M_2}).(1 + c_1(T_{M_B/M_2})t).\]
This shows $\pi_1^*c_5(T_{M_1}).c_1(T_{M_B/M_1}) = \pi_2^*c_5(T_{M_2}).c_1(T_{M_B/M_2})$. Since $c_5(T_{M_i})$ are points, the pull--backs $\pi_i^*c_5(T_{M_i})$ are fibers $F_{i,j}$ of $\pi_i$. Hence
 \[\pi_1^*c_5(T_{M_1}).c_1(T_{M_B/M_1}) = \big(\sum\nolimits_{j=1}^{\deg c_5(T_{M_1})}F_{1,j}\big).(-3{\LB}_1 + 2 {\LB}_2) = 2 \deg c_5(T_{M_1}),\]
since $F_{1,j}.{\LB}_1 = 0$ and $F_{1,j}.{\LB}_2 = 1$ for all $j$. Analogously, $\pi_2^*c_5(T_{M_2}).c_1(T_{M_B/M_2}) = 2 \deg c_5(T_{M_2})$, implying $c_5(T_{M_1}) = c_5(T_{M_2})$, viewed as natural numbers.

It remains hence to compute $c_5(Q)$, where $Q \subset \PN^6$ denotes the $5$--dimensional quadric. From the tangent sequence we get
 \[c_t(T_Q).(1 + 2c_1(\O_Q(1))t) = c_t(T_{\PN^6}|_Q) = (1 + c_1(\O_Q(1))t)^7,\]
where $\O_Q(1) = \O_{\PN^6}(1)|_Q$. Successively we obtain
 \[c_5(T_Q) = \sum\nolimits_{i=0}^5 {7 \choose i} 2^{5-i}c_1(\O_Q(1))^5 = 3c_1(\O_Q(1))^5.\]
Now $c_1(\O_Q(1))^5 = (c_1(\O_{\PN^6}(1)|_Q))^5 = c_1(\O_{\PN^6}(1))^5.Q = c_1(\O_{\PN^6}(1))^5. c_1(\O_{\PN^6}(2)) = 2$, completing the proof.
\end{proof}


\newpage
\section*{Appendix}

The following classification of Fano threefolds with Picard number one is due to Iskovskikh and Mukai, and can be found in \cite{Isk} and \cite{Mukai}, respectively.

\vspace{0.2cm}
\begin{center}
\renewcommand{\arraystretch}{1.3}
\begin{tabular}{*{4}{c|}p{8.5cm}|}
   No.  & $r$ & $H^3$ & $g$ & Description\\\hline
   $1$  & $4$ & $1$ & $33$ & $\PN^3$\\\hline
   $2$  & $3$ & $2$ & $28$ & $Q_3 \subset \PN^4$ the quadric\\\hline
   $3$  & $2$ & $1$ & $5$  & Hypersurface of degree $6$ in $\PN(1,1,1,2,3)$\\
   $4$  & $2$ & $2$ & $9$  & Hypersurface of degree $4$ in $\PN(1,1,1,1,2)$\\
   $5$  & $2$ & $3$ & $13$ & $V_{2,3} \subset \PN^4$ a cubic\\
   $6$  & $2$ & $4$ & $17$ & $V_{2,4} \subset \PN^5$ the intersection of $2$ quadrics\\
   $7$  & $2$ & $5$ & $21$ & $V_{2,5} \subset \PN^6$ the intersection of the Grassmannian ${\rm Gr}(2,5) \subset \PN^9$ by a $\PN^6$\\\hline
   $8$  & $1$ & $2$ & $2$  & Hypersurface of degree $6$ in $\PN(1,1,1,1,3)$\\
   $9$  & $1$ & $4$ & $3$  & $V_4 \subset \PN^4$ a quartic\\
   $10$ & $1$ & $4$ & $3$  & Complete intersection of a quadratic cone and a hypersurface of degree $4$ in $\PN(1,1,1,1,1,2)$\\ 
   $11$ & $1$ & $6$ & $4$ & $V_6 \subset \PN^5$ the intersection of a quadric and a cubic\\
   $12$ & $1$ & $8$ & $5$ & $V_8 \subset \PN^6$ the intersection of three quadrics\\
   $13$ & $1$ & $10$ & $6$ & $V_{10} \subset \PN^7$ the intersection of the Grassmannian ${\rm Gr}(2,5) \subset \PN^9$ by a $\PN^7$\\
   $14$ & $1$ & $12$ & $7$ & $V_{12} \subset \PN^8$ the intersection of the Hermitian symmetric space $M = G/P \subset \PN^{15}$ of type DIII by a $\PN^8$\\
   $15$ & $1$ & $14$ & $8$ & $V_{14} \subset \PN^9$ the intersection of the Grassmannian ${\rm Gr}(2,6) \subset \PN^{14}$ by a $\PN^9$\\
   $16$ & $1$ & $16$ & $9$ & $V_{16} \subset \PN^{10}$ is the intersection of the Hermitian symmetric space $M = G/P \subset \PN^{19}$ of type CI by a $\PN^{10}$\\
   $17$ & $1$ & $18$ & $10$ & $V_{18} \subset \PN^{11}$ is the intersection the $5$--dimensional rational homogeneous contact manifold $G_2/P \subset \PN^{13}$ by a $\PN^{11}$\\
   $18$ & $1$ & $22$ & $12$ & $V_{22} \subset \PN^{13}$ is the zero locus of three sections of the rank $3$ vector bundle $\bigwedge^2{\mathcal Q}$, where ${\mathcal Q}$ is the universal quotient bundle on ${\rm Gr}(7,3)$\\\hline
 \end{tabular}
 \renewcommand{\arraystretch}{1}
\end{center}
\vspace{0.2cm}


\end{document}